\documentclass[CEJP,DVI]{my_cej}

\usepackage{layout}
\usepackage{amsmath,amssymb}
\usepackage{textcomp}
\usepackage{hyperref}
\usepackage{graphicx,slashbox}
\usepackage{subfigure}

\def\w{\left(\omega_k^\a\right)}
\def\a{\alpha}

\def\LDa{{_aD_t^\a}}

\def\LD0{{_0D_t^{0.5}}}

\title{A discrete time method to the first variation of fractional order variational functionals}

\articletype{Research Article}

\author{Shakoor~Pooseh\email{spooseh@ua.pt},
Ricardo~Almeida\email{ricardo.almeida@ua.pt},
Delfim~F.~M.~Torres\email{delfim@ua.pt}}

\institute{CIDMA --- Center for Research and Development in Mathematics and Applications,\\
Department of Mathematics, University of Aveiro, 3810-193 Aveiro, Portugal}

\abstract{The fact that the first variation of a variational functional must vanish along an extremizer
is the base of most effective solution schemes to solve problems of the calculus of variations.
We generalize the method to variational problems involving fractional order derivatives.
First order splines are used as variations, for which fractional derivatives are known.
The Gr\"{u}nwald--Letnikov definition of fractional derivative is used,
because of its intrinsic discrete nature that leads to straightforward approximations.}

\keywords{fractional calculus \*\ numerical approximation \*\ discrete time \*\ fractional calculus of variations}

\pacs{02.30.Xx, 02.60.-x, 45.10.Db, 45.10.Hj}

\begin{document}

\maketitle


\section{Introduction}

Fractional calculus of variations and fractional optimal control deal with dynamic optimization problems
that may depend on some fractional differential equations or involve some fractional operators,
e.g., Riemann--Liouville fractional integrals and derivatives and Caputo fractional derivatives \cite{book:FCV}.
The most common approach to solve such problems consists in using indirect methods: instead of working
directly with the functional, one seeks necessary or sufficient conditions that extremizers must fulfill.
These conditions are often expressed by Euler--Lagrange type equations \cite{Agrawal0,Agrawal1,book:FCV}.
By solving such fractional differential equations, one obtains possible solutions
for the optimization problem. The main disadvantage that rises from such approach,
as it is frequently noted, is that solving analytically those equations is, in many situations,
impossible \cite{MR2727133}. To overcome the problem, numerical methods are applied
\cite{Baleanu0,book:Baleanu,MyID:225}. The most common procedure is to apply numerical techniques
to the fractional differential equations and find an approximation to the extremals.
For example, one can approximate the fractional differential equations
by an ordinary system and apply known techniques from numerical analysis to solve such systems.
In this paper we follow a different approach by neglecting the necessary optimality conditions,
and working directly with the functional.

If $x$ is an extremizer of a given variational functional $J$, the first variation of $J$
evaluated at $x$, along any variation $\eta$, must vanish.
With a discretization on time, we write a system
of equations, being able to solve the initial problem. The advantage of this procedure,
when compared with a discretization method directly on the functional,
is that for different variation functions, one deduce different systems.
We deal with the Riemann--Liouville fractional derivative of order $\a\in(0,1)$
of an integrable function $x:[a,b]\to\mathbb R$, which is defined
(see, e.g., \cite{Kilbas}) by the integral
$$
\LDa x(t)=\frac{1}{\Gamma(1-\alpha)}\frac{d}{dt}\int_a^t \frac{x(\tau)}{(t-\tau)^{\alpha}} d\tau,
$$
where $\Gamma(\cdot)$ is the Gamma function, that is,
$$
\Gamma(z)=\int_0^\infty e^{-t}t^{z-1}\,dt,
$$
but similar techniques can be applied successfully to functionals depending
on other types of fractional derivatives. For instance, using the relation
between the Riemann--Liouville and the Caputo fractional derivatives,
$$
{_a^CD_t^\alpha}x(t)=\LDa x(t)-\frac{x(a)}{\Gamma(1-\alpha)}(t-a)^{-\alpha},
$$
a similar method can be applied to solve fractional variational problems in the Caputo sense.

A regular discretization $t_i=a+ih$, $i=0,\ldots,n$,
over time, with a fixed parameter $h>0$ sufficiently small, is considered,
and we approximate $\LDa x$ by the following finite sum, which is a first
order approximation for the Riemann--Liouville derivative \cite{Podlubny}:
\begin{equation}
\label{approxDerivative}
\LDa x(t_i)\approx \frac{1}{h^\a} \sum_{k=0}^{i}\w x(t_i-kh),
\end{equation}
where
$$
\w=(-1)^k\binom{\a}{k}=(-1)^k\frac{\Gamma(\a+1)}{\Gamma(\a-k+1)k!}.
$$
Throughout the text, the accuracy of each approximation is measured using the maximum norm.
Suppose that $x_i$ is an approximation for $x(t_i)$ on mesh points $t_i$, $i=0,\ldots,n$.
The error caused by the approximation is indicated by $E$ and is computed by
$$
E=\max_{i}|x(t_i)-x_i|.
$$


\section{A numerical method to the fractional calculus of variations}
\label{sec:2}

The problem under consideration is stated in the following way: find the extremizers of
\begin{equation}
\label{def:funct}
J(x)=\int_a^b L\left(t,x(t),\LDa x(t)\right)\,dt
\end{equation}
subject to given boundary conditions $x(a)=x_a$ and $x(b)=x_b$.
Here, $L:[a,b]\times\mathbb R^2\to\mathbb R$ is such that
$\frac{\partial L}{\partial x}$ and $\frac{\partial L}{\partial \LDa x}$
exist and are continuous for all triplets $(t,x(t),\LDa x(t))$.
If $x$ is a solution to the problem and $\eta:[a,b]\to\mathbb R$
is a variation function, i.e., $\eta(a)=\eta(b)=0$,
then the first variation of $J$ at $x$, with the variation $\eta$,
whatever choice of $\eta$ is taken, must vanish:
\begin{equation}
\label{1variation}
J'(x,\eta)=\int_a^b\left[\frac{\partial L}{\partial x}(t,x(t),\LDa x(t))\eta(t)
+\frac{\partial L}{\partial \LDa x}(t,x(t),\LDa x(t))\LDa \eta(t)\right]\,dt=0.
\end{equation}
Using an integration by parts formula for fractional derivatives and the Dubois--Reymond lemma,
Riewe \cite{Riewe:1997} proved that if $x$ is an extremizer of \eqref{def:funct}, then
$$
\frac{\partial L}{\partial x}(t,x(t),\LDa x(t))
+{_tD^\a_b}\left(\frac{\partial L}{\partial \LDa x}\right)(t,x(t),\LDa x(t))=0
$$
(see also \cite{Agrawal}). This fractional differential equation
is called an Euler--Lagrange equation.
For the state of the art on the subject we refer the reader to the recent book \cite{book:FCV}.
Here, instead of solving such Euler--Lagrange equation, we apply a discretization over time
and solve a system of algebraic equations. The procedure has proven
to be a successful tool for classical variational problems \cite{Gregory1,Gregory2}.

The discretization method is the following. Let $n\in\mathbb N$
be a fixed parameter and $h=\frac{b-a}{n}$. If we define
$t_i=a+ih$, $x_i=x(t_i)$, and $\eta_i=\eta(t_i)$ for $i=0,\ldots,n$,
the integral \eqref{1variation} can be approximated by the sum
$$
J'(x,\eta)\approx h \sum_{i=1}^n\left[\frac{\partial L}{\partial x}(t_i,x(t_i),
{_aD_{t_i}^\a} x(t_i))\eta(t_i)+\frac{\partial L}{\partial \LDa x}(t_i,x(t_i),
{_aD_{t_i}^\a} x(t_i)){_aD_{t_i}^\a} \eta(t_i)\right].
$$
To compute the fractional derivative, we replace it by the sum as in
\eqref{approxDerivative}, and to find an approximation for $x$
on mesh points one must solve the equation
\begin{equation}
\label{VarEquation}
\sum_{i=1}^n\left[\frac{\partial L}{\partial x}\left(t_i,x_i,\frac{1}{h^\a}
\sum_{k=0}^{i}\w x_{i-k}\right)\eta_i
+\frac{\partial L}{\partial \LDa x}\left(t_i,x_i,\frac{1}{h^\a}
\sum_{k=0}^{i}\w x_{i-k}\right)\frac{1}{h^\a} \sum_{k=0}^{i}\w \eta_{i-k}\right]=0.
\end{equation}
For different choices of $\eta$, one obtains different equations. Here we use simple variations.
More precisely, we use first order splines as the set of variation functions:
\begin{equation}
\label{VariationFunction}
\eta_j(t)=\left\{
\begin{array}{ll}
\displaystyle\frac{t-t_{j-1}}{h} & \mbox{ if } t_{j-1}\leq t <t_j,\\
\displaystyle\frac{t_{j+1}-t}{h} & \mbox{ if } t_j \leq t<t_{j+1},\\
0 & \mbox{ otherwise,}
\end{array}\right.
\end{equation}
for $j=1,\ldots,n-1$. We remark that conditions $\eta_j(a)=\eta_j(b)=0$
are fulfilled for all $j$, and that $\eta_j(t_i)=0$ for $i\not=j $ and
$\eta_j(t_j)=1$. The fractional derivative of $\eta_j$ at any point $t_i$
is also computed using approximation \eqref{approxDerivative}:
$$
{_aD_{t_i}^\a}\eta_j(t_i)=\left\{
\begin{array}{ll}
\displaystyle\frac{1}{h^\a}(w_{i-j}^\a) & \mbox{ if } j\leq i,\\
0 & \mbox{ otherwise.}
\end{array}\right.
$$
Using $\eta_j$, $j=1,\ldots,n-1$, and equation \eqref{VarEquation},
we establish the following system of $n-1$ algebraic equations
with $n-1$ unknown variables $x_1,\ldots,x_{n-1}$:
\begin{equation}
\label{system1}
\left\{\begin{array}{l}
\displaystyle\frac{\partial L}{\partial x}\left(t_1,x_1,\frac{1}{h^\a}
\sum_{k=0}^{1}\w x_{1-k}\right)
+\frac{1}{h^\a}\sum_{i=1}^{n}\left[\frac{\partial L}{\partial \LDa x}\left(t_i,
x_i,\frac{1}{h^\a} \sum_{k=0}^{i}\w x_{i-k}\right)(w_{i-1}^\a)\right]=0,\\
\displaystyle\frac{\partial L}{\partial x}\left(t_2,x_2,\frac{1}{h^\a}
\sum_{k=0}^{2}\w x_{2-k}\right)
+\frac{1}{h^\a}\sum_{i=2}^{n}\left[\frac{\partial L}{\partial \LDa x}\left(t_i,
x_i,\frac{1}{h^\a} \sum_{k=0}^{i}\w x_{i-k}\right)(w_{i-2}^\a)\right]=0,\\
\quad \vdots\\
\displaystyle\frac{\partial L}{\partial x}\left(t_{n-1},x_{n-1},
\frac{1}{h^\a} \sum_{k=0}^{n-1}\w x_{n-1-k}\right)
+\frac{1}{h^\a}\sum_{i=n-1}^{n}\left[\frac{\partial L}{\partial \LDa x}\left(t_i,
x_i,\frac{1}{h^\a} \sum_{k=0}^{i}\w x_{i-k}\right)(w_{i-n+1}^\a)\right]=0.
\end{array}\right.
\end{equation}
The solution to \eqref{system1}, if exists, gives an approximation
to the values of the unknown function $x$ on mesh points $t_i$.

We have considered so far the so called fundamental
or basic problem of the fractional calculus of variations \cite{book:FCV}.
However, other types of problems can be solved applying similar techniques.
Let us show how to solve numerically the isoperimetric problem, that is,
when in the initial problem the set of admissible functions must satisfy
some integral constraint that involves a fractional derivative.
We state the fractional isoperimetric problem as follows.

Assume that the set of admissible functions are subject not only
to some prescribed boundary conditions, but to some integral constraint, say
$$
\int_a^b g\left(t,x(t),\LDa x(t)\right)\,dt=K,
$$
for a fixed $K\in\mathbb R$. As usual, we assume that
$g:[a,b]\times\mathbb R^2\to\mathbb R$ is such that
$\frac{\partial g}{\partial x}$ and $\frac{\partial g}{\partial \LDa x}$
exist and are continuous. The common procedure to solve this problem
follows some simple steps: first we consider the auxiliary function
\begin{equation}
\label{eq:aux:f}
F=\lambda_0 L(t,x(t),\LDa x(t))+\lambda g(t,x(t),\LDa x(t)),
\end{equation}
for some constants $\lambda_0$ and $\lambda$ to be determined later. Next,
it can be proven that $F$ satisfies the fractional Euler--Lagrange equation
and that in case the extremizer does not satisfies the Euler--Lagrange associated to $g$,
then we can take $\lambda_0=1$ (cf. \cite{Almeida1}).
In conclusion, the first variation of $F$ evaluated
along an extremal must vanish, and so we obtain a system similar to \eqref{system1},
replacing $L$ by $F$. Also, from the integral constraint, we obtain another
equation derived by discretization that is used to obtain $\lambda$:
$$
h \sum_{i=1}^n g\left(t_i,x_i,\frac{1}{h^\a} \sum_{k=0}^{i}\w x_{i-k}\right)=K.
$$


\section{Illustrative examples}

We show the usefulness of our approximate method
with three problems of the fractional calculus of variations.


\subsection{Basic fractional variational problems}

\begin{example}
\label{Ex1}
Consider the following variational problem: to minimize the functional
$$
J(x)=\int_0^1\left(\LD0 x(t)-\frac{2}{\Gamma(2.5)}t^{1.5}\right)^2\,dt
$$
subject to the boundary conditions $x(0)=0$ and $x(1)=1$.
It is an easy exercise to verify that the solution is the function $x(t)=t^2$.
We apply our method to this problem, for the variation \eqref{VariationFunction}.
The functional $J$ does not depend on $x$ and is quadratic with respect
to the fractional term. Therefore, the first variation is linear.
The resulting algebraic system from \eqref{system1} is also linear and easy to solve:
$$
\left\{\begin{array}{l}
\displaystyle \sum_{i=0}^{n-1}\left(\omega_i^{0.5}\right)^2 x_1
+\sum_{i=1}^{n-1}\left(\omega_i^{0.5}\right)\left(\omega_{i-1}^{0.5}\right) x_2
+\sum_{i=2}^{n-1}\left(\omega_i^{0.5}\right)\left(\omega_{i-2}^{0.5}\right)x_3\\
\quad \displaystyle +\cdots+
\sum_{i=n-2}^{n-1}\left(\omega_i^{0.5}\right)\left(\omega_{i-(n-2)}^{0.5}\right)x_{n-1}
=\frac{2h^{2}}{\Gamma(2.5)}\sum_{i=0}^{n-1}\left(\omega_i^{0.5}\right)(i+1)^{1.5}
-\left(\omega_0^{0.5}\right)\left(\omega_{n-1}^{0.5}\right),\\
\displaystyle \sum_{i=0}^{n-2}\left(\omega_i^{0.5}\right)\left(\omega_{i+1}^{0.5}\right) x_1
+\sum_{i=0}^{n-2}\left(\omega_i^{0.5}\right)^2 x_2
+\sum_{i=1}^{n-2}\left(\omega_i^{0.5}\right)\left(\omega_{i-1}^{0.5}\right) x_3\\
\quad \displaystyle +\cdots +
\sum_{i=n-3}^{n-2}\left(\omega_i^{0.5}\right)\left(\omega_{i-(n-3)}^{0.5}\right)x_{n-1}
=\frac{2h^{2}}{\Gamma(2.5)}\sum_{i=0}^{n-2}\left(\omega_i^{0.5}\right)(i+2)^{1.5}
-\left(\omega_0^{0.5}\right)\left(\omega_{n-2}^{0.5}\right),\\
\qquad \vdots\\
\displaystyle \sum_{i=0}^{1}\left(\omega_i^{0.5}\right)\left(\omega_{i+n-2}^{0.5}\right) x_1
+\sum_{i=0}^{1}\left(\omega_i^{0.5}\right)\left(\omega_{i+n-3}^{0.5}\right) x_2
+\sum_{i=0}^{1}\left(\omega_i^{0.5}\right)\left(\omega_{i+n-4}^{0.5}\right) x_3\\
\quad \displaystyle +\cdots + \sum_{i=0}^{1}\left(\omega_i^{0.5}\right)^2x_{n-1}
=\frac{2h^{2}}{\Gamma(2.5)}\sum_{i=0}^{1}\left(\omega_i^{0.5}\right)(i+n-1)^{1.5}
-\left(\omega_0^{0.5}\right)\left(\omega_{1}^{0.5}\right).
\end{array}\right.
$$
The exact solution together with three numerical approximations,
with different discretization step sizes,
are depicted in Figure~\ref{FigEx1}.
\begin{figure}[!h]
\center
\includegraphics[scale=0.5]{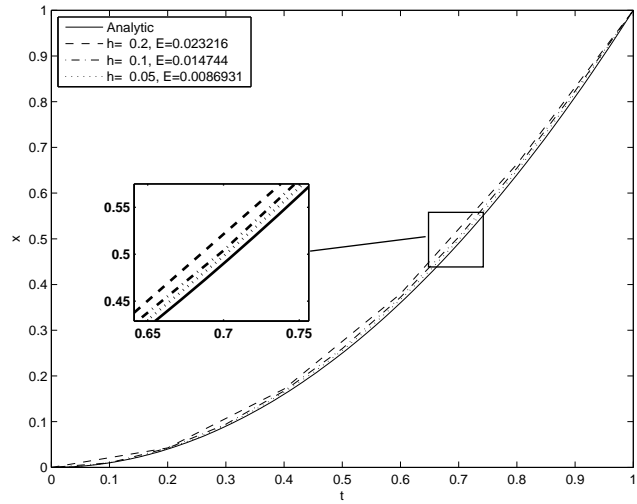}
\caption{Exact solution versus numerical approximations to Example~\ref{Ex1}.}
\label{FigEx1}
\end{figure}
\end{example}

\begin{example}
\label{Ex2}
Find the minimizer of the functional
$$
J(x)=\int_0^1\left(\LD0 x(t)-\frac{16\Gamma(6)}{\Gamma(5.5)}t^{4.5}+
\frac{20\Gamma(4)}{\Gamma(3.5)}t^{2.5}-\frac{5}{\Gamma(1.5)}t^{0.5}\right)^4\,dt,
$$
subject to $x(0)=0$ and $x(1)=1$. The minimum value of this functional is zero and the minimizer is
$$
x(t)=16t^{5}-20t^{3}+5t.
$$
Discretizing the first variation as discussed above, leads to a nonlinear system of algebraic equation.
Its solution, using different step sizes, is depicted in Figure~\ref{FigEx2}.
\begin{figure}[!h]
\center
\includegraphics[scale=0.5]{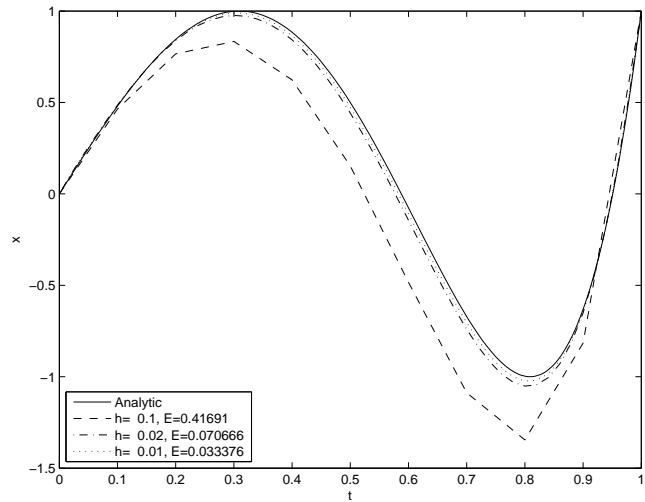}
\caption{Exact solution versus numerical approximations to Example~\ref{Ex2}.}
\label{FigEx2}
\end{figure}
\end{example}


\subsection{An isoperimetric fractional variational problem}
\label{Isoperimetricproblem}

\begin{example}
\label{Ex3}
Let us search the minimizer of
$$
J(x)=\int_0^1\left(t^4 +\left(\LD0 x(t)\right)^2\right)\,dt
$$
subject to the boundary conditions
$$
x(0)=0 \quad \mbox{and}\quad x(1)=\frac{16}{15\Gamma(0.5)}
$$
and the integral constraint
$$
\int_0^1t^2\,\LD0 x(t)\,dt=\frac{1}{5}.
$$
In \cite{Almeida1} it is shown that the solution to this problem is the function
$$
x(t)=\frac{16t^{2.5}}{15\Gamma(0.5)}.
$$
Because $x$ does not satisfy the fractional Euler--Lagrange equation associated
to the integral constraint, one can take $\lambda_0=1$
and the auxiliary function \eqref{eq:aux:f} is
$F=t^4 +\left(\LD0 x(t)\right)^2+\lambda \, t^2\,\LD0 x(t)$.
Now we calculate the first variation of $\int_0^1 F dt$. An extra unknown,
$\lambda$, is present in the new setting, that is obtained by discretizing
the integral constraint, as explained in Section~\ref{sec:2}.
The solutions to the resulting algebraic system,
with different step sizes, are given in Figure~\ref{FigIso}.
\begin{figure}[!h]
\centering
\includegraphics[scale=0.5]{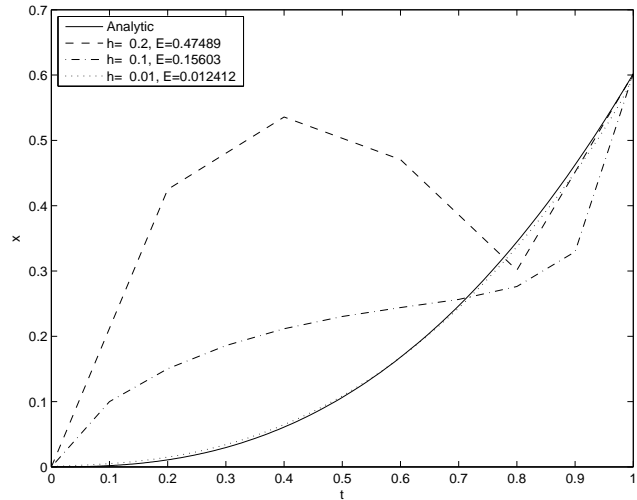}
\caption{Exact versus numerical approximations
to the isoperimetric problem of Example~\ref{Ex3}.}
\label{FigIso}
\end{figure}
\end{example}


\section*{Acknowledgments}

Work supported by {\it FEDER} funds through {\it COMPETE}
-- Operational Programme Factors of Competitiveness
(''Programa Operacional Factores de Competitividade'')
and by Portuguese funds through the {\it Center for Research
and Development in Mathematics and Applications} (University of Aveiro)
and the Portuguese Foundation for Science and Technology
(''FCT -- Funda\c{c}\~{a}o para a Ci\^{e}ncia e a Tecnologia''),
within project PEst-C/MAT/UI4106/2011 with COMPETE
number FCOMP-01-0124-FEDER-022690.
Part of first author's Ph.D.,
which is carried out at the University of Aveiro under the
\emph{Doctoral Program in Mathematics and Applications} (PDMA)
of Universities of Aveiro and Minho, supported
by the Ph.D. fellowship SFRH/BD/33761/2009.



\end{document}